\documentclass[fleqn]{mat01}
\usepackage{times,mathtimy,amssymb,latexsym}%,krepsf,rangecite}
\begin{document}

\setcounter{page}{83}
\firstpage{83}

\def\d{\mbox{\rm d}}
\def\e{\mbox{\rm e}}

\newtheorem{defini}{\rm DEFINITION}[ssection]
\newtheorem{corol}{\rm COROLLARY}[ssection]
\newtheorem{exa}{\it Example}[ssection]

\renewcommand{\theequation}{\arabic{section}.\arabic{equation}}
\def \qed {{\nopagebreak}\hfill \vrule height6pt width 6pt depth 0pt}

\newtheorem{assump}{{\it Assumption A$\!\!$}}
\newcommand{\bhyp}{\begin{assump}}
\newcommand{\ehyp}{\end{assump}}
\newcommand{\raro}{\rightarrow}
\newcommand {\Raro}{\Rightarrow}
\newcommand{\lra}{\Longrightarrow}

\newcommand{\lng}{\langle}
\newcommand{\rng}{\rangle}

\newcommand{\clh}{{\cal H}}
\newcommand{\clm}{{\cal M}}
\newcommand{\clg}{{\cal G}}
\newcommand{\cls}{{\cal S}}
\newcommand{\clc}{{\cal C}}
\newcommand{\fal}{\forall}
\newcommand{\cld}{{\cal D}}
\newcommand{\IP}{\mathbb{P}}
\newcommand{\sbs}{\subset}
\newcommand{\clb}{{\cal B}}
\newcommand{\clp}{{\cal P}}
\newcommand{\cll}{{\cal L}}
\newcommand{\clf}{{\cal F}}
\newcommand{\IE}{\mathbb{E}} %%{I\!\!E}
\newcommand{\ity}{\infty}
\newcommand{\clw}{{\cal W}}

\def\wtit{\widetilde}
\newcommand{\IQ}{\,\mathbb{Q}\ } %%{I\!\!\!Q}
\newcommand{\IR}{\mathbb{R}}

\title{On characterisation of Markov processes via martingale problems}

\markboth{Abhay G Bhatt, Rajeeva L Karandikar and B V Rao}{Characterisation of Markov processes}

\author{ABHAY G BHATT, RAJEEVA L KARANDIKAR and B V RAO$^{*}$}

\address{Indian Statistical Institute, 7, SJS Sansanwal Marg, New Delhi~110~016, India\\
\noindent $^{*}$Indian Statistical Institute, 203, B.T.~Road, Kolkata~700~108, India\\
\noindent E-mail: abhay@isid.ac.in; rlk@isid.ac.in;
bvrao@isical.ac.in}

\volume{116}

\mon{February}

\parts{1}

\pubyear{2006}

\Date{MS received 31 March 2005; revised 28 September 2005}

\begin{abstract}
It is well-known that well-posedness of a martingale problem in the
class of continuous (or r.c.l.l.) solutions enables one to construct the
associated transition probability functions. We extend this result to
the case when the martingale problem is well-posed in the class of
solutions which are continuous in probability. This extension is used to
improve on a criterion for a probability measure to be invariant for the
semigroup associated with the Markov process. We also give examples of
martingale problems that are well-posed in the class of solutions which
are continuous in probability but for which no r.c.l.l. solution exists.
\end{abstract}


\maketitle

\section{Introduction}

The seminal paper on multi-dimensional diffusions by Stroock and
Varadhan~\cite{SVp} introduced martinagle problems as a way of
construction and study of Markov processes. Since then, this
approach has been used successfully in several contexts such as
interacting particle systems, Markov processes associated with
Boltzmann equation, nonlinear filtering theory, controlled Markov
processes, branching processes etc. A good account of the
`\hbox{theory} of martingale problems' is given in the book by
Ethier and Kurtz~\cite{EK}. To construct a Markov process, the
martingale problem approach allows one to construct the process
for each initial condition separately and a general result gives
the measurability of the associated transition probability
function. To proceed, we give the basic definitions here.

Given an operator $A$ with domain $\cld(A)\subseteq C_b(E)$ and range
subset of $C_b(E)$ (where $E$ is a complete separable metric space), a
process $X^x_t$ adapted to a filtration $({\cal F}_t)$ is said to be a
solution to the $(A,\delta_x)$ martingale problem if for all
$f\in\cld(A)$,
\begin{equation}\label{r1}
f(X^x_t)-\int_0^tAf(X^x_u)\d u \mbox{ is a $({\cal F}_t)$-martingale}
\end{equation}
and
\begin{equation}\label{r2}
\IP(X^x_0=x)=1.
\end{equation}
The martingale problem for $A$ is said to be well-posed in the class of
r.c.l.l. solutions if for all $x$ there exists a r.c.l.l. process
$(X^x_t)$ satisfying (\ref{r1}) and (\ref{r2}) and further for two such
processes satisfying (\ref{r1}) and (\ref{r2}) (defined possibly on
different probability spaces), the finite dimensional distributions are
the same. Well-posedness in the class of continuous solutions or
measurable solutions is similarly defined. A well-known result, which
has its origins in the work of Stroock and Varadhan~\cite{SV} says that if the
martingale problem for $A$ is well-posed in the class of r.c.l.l.
solutions (or well-posed in the class of continuous solutions), then
(assuming that $A$, $\cld(A)$ satisfy some mild conditions) it follows
that $p_t(x,\cdot)$ defined by
\begin{equation} \label{r3}
p_t(x,A)=\IP(X^x_t\in A)
\end{equation}
is a transition probability function and any solution is a Markov
process with $p_t$ as its transition probability function (see
e.g. Theorems IV.4.2 and IV.4.6 of \cite{EK}). This in
turn gives us the associated semi-group $(T_t)$ and its generator $L$.
The generator $L$ happens to be an extension of the operator $A$ and
thus $A$ contains all the `relevant information' about $L$ as well as
about $X$.

We extend this result and show that if the martingale problem is
well-posed in the class of solutions that are continuous in
probability, then (under suitable conditions on $A,\cld(A)$) the
function $p_t$ defined by (\ref{r3}) is measurable.

In order to achieve our aim we give a Borel structure to the set of
distributions of processes that are continuous in probability. Once we
have done this, we can deduce that well-posedness in the class of
solutions which are continuous in probability implies measurability of
the associated transition probability function.

In \S4, we give criterion for a measure to be invariant for the
semigroup generated by a well-posed martingale problem. This is an
improvement on several results on this theme (see \cite{BB,RLKAB1,RLKAB2,RLKAB3,Ech,KS1,KS2}). In the last section, we give examples of
operators (and their domains) satisfying the conditions of \S3,
and such that the corresponding martingale problems are well-posed in
the class of solutions that are continuous in probability but for which
no r.c.l.l. solution exists.

\section{Preliminaries}

We will denote by $(E, d)$ a complete, separable metric space. $A$
will denote an operator with domain $D(A) \sbs C_b(E)$, the space
of real-valued bounded continuous functions on $E$ and with range
contained in $M(E)$, the class of all real-valued Borel measurable
functions on $E$. Let $B(E)$ denote the class of all bounded Borel
measurable functions. For $C\sbs B(E)$, we define the $bp$-closure
of $C$ to be the smallest subset of $B(E)$ containing $C$ which is
closed under bounded pointwise convergence of sequences of
functions. $\clb (E)$ will denote the Borel $\sigma$-field on $E$,
$\clp (E)$ will denote the space of probability measures on $E$.
For a random variable $Z$ taking values in $E$, $\cll(Z)$ will
denote the law of $Z$- i.e. the probability measure $\IP\circ
Z^{-1}$, if $Z$ is defined on $(\Omega , \clf , \IP)$. For a
measurable Process$(X_t)$ defined on $(\Omega , \clf , \IP)$, let
\begin{equation*}
^{*}\clf^X_t= \sigma\left\lbrace X_u, \; \int_0^u h(X_s)\d
s\hbox{:}\;u\leq t\;,\;h\in C_b(E)\right\rbrace.
\end{equation*}

Throughout this article, we will assume the following:
\bhyp\label{h1}
{\rm There exists a $[0,\infty)$-valued measurable function $\Phi$ on $E$ such that
\setcounter{equation}{0}
\begin{equation}
|Af(x)|\leq C_f\Phi(x)\;\;\;\fal x\in E,\; f\in
\cld(A).\label{rk1}
\end{equation}}
\ehyp

\begin{defini}$\left.\right.$\vspace{.5pc}

\noindent {\rm An $E$-valued process $(X_t)_{0\leq t<\infty}$
defined on some probability space $(\Omega , \clf , \IP)$ is said
to be a solution to the martingale problem for $(A, \mu )$ if

\begin{enumerate}
\leftskip .4pc
\renewcommand{\labelenumi}{(\roman{enumi})}
\item  $X$ is a measurable process with $\cll (X_0) = \mu $,
\item $\IE_\IP\left[\int_0^T \Phi (X_s)\d s\right] < \ity $ for all $T<\infty$,
and
\item for every $f \in \cld(A) $,
\begin{equation}\label{mft}
\hskip -1.25pc M^f_t=f(X_t ) - \int^t_0 Af(X_s)\ \d s
\end{equation}
is a $({}^*\clf^X_t)$-martingale. $X$ will be called a solution to the $A$
martingale problem if it is a solution to the $(A,\mu)$ martingale problem
for some $\mu$.
\end{enumerate}

Let $\clw$ be a class of $E$-valued processes. For example, we could
consider $\clw$ to be the class of $E$-valued processes with r.c.l.l.
paths or $\clw$ can be the class of solutions that are continuous in
probability.}
\end{defini}

\begin{defini}$\left.\right.$\vspace{.5pc}

\noindent {\rm The martingale problem for $A$ is said to be
well-posed in the class $\clw$ if for all $x\in E$, there exists a
solution $X^x \in \clw$ to the $(A,\delta_x)$ martingale problem
and if $Y\in\clw$ is any other solution to the $(A,\delta_x)$
martingale problem, then the finite dimensional distributions of
$X^x$ and $Y$ are the same.

We begin with some observations on solutions to the $A$-martingale problem.}
\end{defini}

\begin{theorem}[\!]\label{NT1}
Let $X$ {\rm (}defined on some $(\Omega , \clf , \IP)${\rm )} be a solution of the
martingale problem for $A$. Suppose that $\cld(A)$ is a determining
class and further that
\begin{equation}\label{contlaw}
t\longrightarrow \cll(X_t)\mbox{ is continuous.}
\end{equation}
Then $t\longrightarrow X_t$ is continuous in probability.
\end{theorem}

\begin{proof}
Let $f\in\cld(A)$. The assumption (\ref{contlaw}) alongwith the
fact that the martingale $M^f$ (see eq.~\eqref{mft}) has a
r.c.l.l. modification implies that $M^f$ is continuous in
probability. This in turn implies that the mapping
$t\longrightarrow f(X_t)$ is continuous in probability. As a
consequence, for $f,g\in\cld(A)$,
\begin{equation}\label{r13}
(s,t)\longrightarrow \IE_\IP[f(X_s)g(X_t)] \mbox{ is continuous}.
\end{equation}
Let $s_k\raro s$. The assumption (\ref{contlaw}) implies that the family
of distributions $\{\cll(X_{s_k})\}$ is tight and so the family of
distributions (on $E\times E$)
\begin{equation}\label{r14}
\{\cll(X_{s_k},X_{s})\hbox{:}\ k\geq 1\} \mbox{ is tight}.
\end{equation}
Since the class of functions $(x,y)\longrightarrow f(x)g(y)$,
$f,g\in\cld(A)$ constitutes a determining class, (\ref{r13})
and (\ref{r14}) together imply that
 \begin{equation}\label{r15}
 \cll(X_{s_k},X_{s})\raro \cll(X_{s},X_{s}).
 \end{equation}
Now for any $\epsilon>0$, $\IP(d(X_{s},X_{s})\geq \epsilon)=0$.
Thus in view of (\ref{r15})
\begin{equation*}
\limsup_{k\raro\infty} \IP(d(X_{s_k},X_{s})\geq \epsilon)\leq 0
\end{equation*}
i.e.,
\begin{equation*}
\IP(d(X_{s_k},X_{s})\geq \epsilon)\raro 0.
\end{equation*}
This completes the proof.
\qed
\end{proof}

\begin{remark} \label{remark1}
{\rm The proof given above contains the proof of the following: if for a
process $Y$, the mapping $(s,t)\longrightarrow \cll(Y_s,Y_t)$ is
continuous, then $Y$ is continuous in probability.}
\end{remark}
\begin{remark}
{\rm The assumption (\ref{contlaw}) can be replaced by
\begin{equation*}
\{X_t\hbox{:}\ 0\leq t \leq T\} \mbox{ is tight } \forall
T<\infty.
\end{equation*}}
\end{remark}

\setcounter{corol}{3}
\begin{corol}\label{NC2}$\left.\right.$\vspace{.5pc}

\noindent Let $X$ {\rm (}defined on some $(\Omega , \clf , \IP)${\rm )} be a solution of the
martingale problem for $A$. Suppose that the domain $\cld(A)$ of $A$ is
a convergence determining class on $E$. Then the process $X$ is
continuous in probability.
\end{corol}

\begin{proof}
Since $f(X_t)-\int_0^tAf(X_s) \d s$ is a martingale for
$f\in\cld(A)$, it follows that the mapping $t\longrightarrow
\IE_\IP[f(X_t)]$ is continuous. Since $\cld(A)$ is a convergence
determining class, this implies continuity of the mapping
\begin{equation*}
t\longrightarrow \cll(X_t).
\end{equation*}
Thus, by Theorem \ref{NT1}, $t\longrightarrow X_t$ is continuous in probability.
\qed
\end{proof}

\section{Main result}
\setcounter{equation}{0}

We have seen in the previous section that under suitable conditions, all
solutions to a martingale problem are continous in probability. Thus we
now construct a Borel structure on the class of distributions of such
processes.

For $m\geq 1$, $E^m$ with the product topology is again a complete
separable metric space. Let $\clp(E^m)$ be equipped with the
topology of weak convergence. Let
$\clc_m=C([0,\infty)^m,\clp(E^m))$ be equipped with the topology
of uniform convergence on compact subsets. Then $\clc_m$ is a
complete separable metric space. Let $\cls_m$ be the set of
$\mu^m=\mu^m(t_1,t_2,\dots,t_m)\in\clc_m$ satisfying
\begin{align}
&\int (\pi f)(x_1,x_2,\dots,x_m) \mu^m(t_{\pi 1},t_{\pi 2},\dots,t_{\pi m})(\d x_1,\d x_2,\dots,\d x_m)\nonumber\\[.5pc]
&\quad\,=\int f(x_1,x_2,\dots,x_m) \mu^m(t_{ 1},t_{ 2},\dots,t_{ m})(\d x_1,\d x_2,\dots,\d x_m)
\end{align}
for all permutations $\pi$ of $\{1,2,\dots ,m\}$, for all $f\in
C_b(E^m)$ where $\pi f$ is defined by
\begin{equation*}
\pi f(x_1,x_2,\dots,x_m)= f(x_{\pi 1},x_{\pi 2},\dots,x_{\pi m}).
\end{equation*}

It is easy to see that $\cls_m$ is a closed subset of $\clc_m$ and hence
$\cls_m$ is a complete separable metric space. Let
$\cls_\infty=\Pi_{m=1}^{\infty}\cls_m$. Under the product topology,
$\cls_\infty$ is also a complete separable metric space. Elements of
$\cls_\infty$ will be denoted by $\mu=(\mu^1,\mu^2,\dots )$ with
$\mu^k\in\cls_k$. Let $D$ denote the diagonal in $E^2$,
\begin{equation*}
D=\{(x,x)\hbox{:}\ x\in E\}
\end{equation*}
and let
\begin{equation*}
\clh=\{\mu^2\in\cls_2\hbox{:}\ \mu^2(t,t)(D)=1\quad \forall t \in
[0,\infty)\}.
\end{equation*}
Since $D$ is closed in $E^2$ and $\mu^2\in\cls_2$ is continuous, it
follows that $\clh$ is a closed subset of $\cls_2$. Let
\begin{align*}
\cls^*=\{\mu\in\cls_\infty\hbox{:}\ \mu^m(t_1,\dots ,t_m)\circ
(h_m)^{-1}=\mu^{m-1}(t_1,\dots ,t_{m-1}),\,\forall m>1\},
\end{align*}
where $h_m\hbox{:}\ E^m\longrightarrow E^{m-1}$ is the projection
map defined by
\begin{equation*}
h_m(x_1,x_2,\dots,x_m)=(x_1,x_2,\dots,x_{m-1}).
\end{equation*}
Let
\begin{equation*}
\tilde{\cls}=\{\mu\in\cls^*\hbox{:}\ \mu^2\in\clh\}.
\end{equation*}
Then clearly $\tilde{\cls}$ is also a complete separable metric
space since it is a closed subspace of $\cls_\infty$. Every
element of $\tilde{\cls}$ is a {\it consistent\/} family of finite
dimensional distributions and hence by the Kolmogorov consistency
theorem, given $\mu=(\mu^1,\mu^2,\dots )\in\tilde{\cls}$, there
exists a probability space $(\Omega^*,\clf^*,P^*)$ and a
stochastic process $(X_t)$ on it such that for all $m\geq 1$,
\begin{equation}\label{r21}
\cll(X_{t_1},X_{t_2},\dots ,X_{t_m})=\mu^m(t_1,t_2,\dots ,t_m).
\end{equation}
In view of Remark \ref{remark1} and the fact that $\mu^2\in\clh$,
the process $X$ is continuous in probability. Conversely, given a
$E$-valued process $X$ that is continuous in probability, $\mu^m$
defined by (\ref{r21}) belongs to $\cls_m$, $\mu^2$ belongs to
$\clh$ and clearly $\{\mu^1,\mu^2,\dots \}$ is a consistent family
and hence $\mu=(\mu^1,\mu^2,\dots ) \in\tilde{\cls}$. Thus,
$\tilde{\cls}$ can be identified with the class of distributions
of $E$-valued processes that are continuous in probability.

Having given a topological structure to the class of (distributions of)
processes that are continuous in probability, we now identify the class
of (distributions of) solutions to the martingale problem for $A$ and
show that it is a Borel set. As in the corresponding result on solutions
with r.c.l.l. paths \cite{EK}, we assume that
$A,\cld(A)$ satisfy the following:

\bhyp \label{h2} {\rm There exists a countable set $\{
f_n\hbox{:}\ n\geq 1 \} \sbs \cld(A)$ such that
\begin{equation*}
bp-{\rm closure}\{ (f_n, \Phi ^{-1}  Af_n)\hbox{:}\ n \geq 1\}
\supset \{ (f, \Phi ^{-1}Af)\hbox{:}\ f \in \cld(A) \}.
\end{equation*}}
\ehyp
Let $X$ be a process that is continuous in probability (on some
$(\Omega,\clf,\IP))$. Since every such process admits a measurable
modification \cite{DM}, we assume that $X$ is
measurable. Let $\clg$ be a countable $bp$-dense subset of $C_b(E)$. Then
$X$ is a solution to the $A$ martingale problem if and only if
\begin{equation*}
\IE_\IP\Bigl[\int_0^N\phi(X_u)\d u\Bigr]<\infty\quad \forall N\geq
1
\end{equation*}
and
\begin{align*}
\IE_\IP\left[g_1(X_{s_1})\dots g_k(X_{s_k})\left(f_m(X_t)-f_m(X_s)-\int_s^t(Af_m)(X_u)\d u\right)\right]=0
\end{align*}
for all $s_1, s_2,\dots,s_k, s,t$ rationals with $s_i\leq s\leq
t$, $g_i\in\clg$, $1\leq i\leq k$, $k\geq 1,m\geq 1$, where
$\{f_j\hbox{:}\ j\geq 1\}$ are as in Assumption A2. Thus, a
measurable process $X$ is a solution to the $A$ martingale problem
if and only if its finite dimensional distributions
$\mu=(\mu^1,\mu^2,\dots )$ defined by (\ref{r21}) belong to
$\clm\sbs\tilde{\cls}$ defined as follows: $\clm$ is the set of
$\mu=(\mu^1,\mu^2,\dots )\in\tilde{\cls}$ satisfying
\begin{equation}\label{r25}
\int_0^N\lng \mu^1(s),\Phi\rng \d s<\infty\quad \forall N\geq 1
\end{equation}
(here, $\lng F,\Gamma\rng$ denotes $\int F\d\Gamma$) and
\begin{align}\label{r26}
&\lng\mu^{k+1}(s_1,s_2,\dots ,s_k,t), G\otimes f_m\rng - \lng\mu^{k+1}(s_1,s_2,\dots ,s_k,s), G\otimes f_m\rng\nonumber\\[.3pc]
&\quad\,=\int_s^t\lng\mu^{k+1}(s_1,s_2,\dots ,s_k,u), G\otimes Af_m\rng \d u
\end{align}
for all $s_1, s_2,\dots s_k, s,t$ rationals with $s_i\leq s\leq
t$, $g_i\in\clg$, $1\leq i\leq k$, $k\geq 1,m\geq 1$, where
$\{f_j\hbox{:}\ j\geq 1\}$ are as in Assumption A2 and
\begin{equation*}
G\otimes f_m(x_1,x_2,\dots ,x_k,z)=g_1(x_1)g_2(x_2)\dots g_k(x_k)f_m(z).
\end{equation*}
Since $\clm$ is defined via countably many conditions with each
condition in turn involving measurable functions of
$\mu=(\mu^1,\mu^2,\dots )$, it follows that $\clm$ is a Borel
subset of $\tilde{\cls}$. Moreover, given $\mu=(\mu^1,\mu^2,\dots
)\in\clm$, as noted above there exists a process $X$ such that its
finite dimensional distributions are those given by
$\mu=(\mu^1,\mu^2,\dots )$. Further, $\mu^2\in\clh$ and Remark 2.2
implies that this process is continuous in probability and can be
assumed to be measurable. It follows that $X$ is a solution to the
$A$ martingale problem. We have thus proved the following.

\setcounter{defin}{0}
\begin{theorem}[\!]
Suppose that $A, \cld(A)$ satisfy Assumptions A{\rm 1} and A{\rm 2}.
Then $\mu=(\mu^1,\mu^2,\dots) \in \clm$ if and only if there exists a
process $X$ that is {\rm (i)} continuous in probability{\rm , (ii)} the finite
dimensional distributions of $X$ are given by $\mu=(\mu^1,\mu^2,\dots
)$ and {\rm (iii)} $X$ is a solution to the martingale problem for $A$.
\end{theorem}

We are now ready to prove the measurability of $p_t$ when the martingale
problem for $A$ is well-posed. We introduce the following:

\bhyp\label{h3}
{\rm The martingale problem for $(A,\delta_x)$ is well-posed in the class of
solutions that are continuous in probability for each $x\in E$.}
\ehyp

\begin{theorem}[\!]\label{mainT}

Suppose that $A,\cld(A)$ satisfy A{\rm 1}, A{\rm 2} and A{\rm 3}.
Let $X^x$ denote a solution that is continuous in probability to the
$(A,\delta_x)$ martingale problem. Let $p_t(x,B),\,t\in
[0,\infty),\,x\in E ,\,B\in\clb(E)$ be defined by
\begin{equation}\label{r31}
p_t(x,B)=P(X^x_t\in B).
\end{equation}
Then for all $t\in [0,\infty),\,B\in\clb(E)${\rm ,} $x\longrightarrow
p_t(x,B)$ is Borel measurable.
\end{theorem}

\begin{proof}
Note that $F=\{\delta_x\hbox{:}\ x\in E\}$ is a Borel measurable
subset of $\clp (E)$ (indeed it is a closed subset) and the
function $\theta(\delta_x)=x$ is a Borel measurable function on it
(again this is a continuous function). Let $\psi_t\hbox{:}\
\clm\longrightarrow \clp (E)$ for $0\leq t<\infty$ be defined by
\begin{equation*}
\psi_t(\mu)=\mu^1(t),\;\;\,\;\;\,\mu=(\mu^1,\mu^2,\dots ) \in\clm.
\end{equation*}
The functions $\psi_t $ are continuous and hence measurable. Let
$\clm_0=(\psi_0)^{-1}(F)$. It follows that $\clm_0$ is a Borel
subset of $\tilde{\cls}$. Also, $\Psi=\theta\circ\psi_0$ is a
measurable function from $\clm_0$ into $E$.

In view of the Assumption A3, for a given $x\in E$, $\clm$ has exactly
one element $\mu=(\mu^1,\mu^2,\dots ) $ such that
\begin{equation*}
\mu^1(0)=\delta_x
\end{equation*}
and hence the function $\Psi$ is one-to-one. Hence by Kurtowski's
theorem (see e.g. Corollary I.3.3 of \cite{KRP}) the
function is bimeasurable, or it has a measurable inverse. Let us note
that $\Psi^{-1}(x)$ denotes the finite dimensional distributions of
$X^x$- the (unique in law) solution to $(A,\delta_x)$ martingale problem
which is continuous in probability. The required conclusion follows by
noting that\vspace{.5pc}

\hfill $p_t(x,B)=\psi_t(\Psi^{-1}(x))(B)$.\hfill
\qed
\end{proof}

\bhyp\label{h4}
{\rm $\cld(A)$ is convergence determining.}
\ehyp

\bhyp\label{h5}
{\rm The $(A,\delta_x)$ martingale problem is well-posed in the class of measurable processes for all $x\in E$.}
\ehyp

\begin{remark}
{\rm Let us note that Assumptions A4 and A5 imply Assumption A3. This is
because Assumption A4 implies that every solution to the $A$ martingale
problem is continuous in probability. Thus the conclusion of the above
theorem remains valid with the same proof if instead we assume that
$A,\cld(A)$ satisfy \mbox{A1}, \mbox {A2} , \mbox{A4} and \mbox{A5}.}
\end{remark}

\begin{remark} \label{mux}
{\rm Assume that Assumptions A1, A2 and A3 are true. Denote by
\begin{equation*}
\mu_x=(\mu_x^1,\mu_x^2,\dots )
\end{equation*}
the finite dimensional distributions of the (unique in law) solution to
the $(A,\delta_x)$ martingale problem that is continuous in probability.
We have seen in the proof above that
\begin{equation*}
x\longrightarrow \mu_x(=\Psi^{-1}(x))
\end{equation*}
is Borel measurable and hence for all $t_1,t_2,\dots ,t_m$, $m\geq
1$
\begin{equation}\label{r33}
x\longrightarrow \mu_x^m(t_1,t_2,\dots ,t_m) \mbox{ is Borel measurable}.
\end{equation}}
\end{remark}

The next step is to prove that $\{T_t\hbox{:}\ t\geq 0\}$ defined
by
\begin{equation}\label{r41}
T_tf(x)=\int f(y)p_t(x,\d y)=\int f(y)\mu_x^1(t)(\d y)
\end{equation}
is a semigroup on the class of bounded Borel measurable functions $f$ on
$E$. For this, we need to consider the martingale problem with
non-degenerate initial distributions. Note that well-posedness for
degenerate initials in the class of all solutions may not imply well-posedness
for all initials. To proceed further, let us introduce the
following notation:
\begin{equation}\label{r42}
\Phi^*_N(x)=\int_0^N\lng\mu_x^1(s),\Phi\rng \d s.
\end{equation}
Then in view of Remark \ref{mux}, it follows that $\Phi^*_N$ is a $[0,\infty)$-valued measurable function. The next lemma shows that the existence of solution
to the martingale problem holds for a large class
of initial distributions.

Let $\clp_\Phi$ be the set of all measures $\lambda\in\clp(E)$ such that
\begin{equation}\label{r43}
\lng\Phi^*_N,\lambda\rng<\infty,\quad \forall N\geq 1.
\end{equation}

\begin{lemma}\label{L7}
Suppose that $A,\cld(A)$ satisfy Assumptions A{\rm 1,} A{\rm 2} and
A{\rm 3}. Let $\lambda\in\clp_\Phi$. Then $\nu=(\nu^1,\nu^2,\dots )$
defined by
\begin{equation}\label{r44}
\lng\nu^m(t_1,t_2,\dots ,t_m),g\rng=\int \lng\mu_x^m(t_1,t_2,\dots ,t_m),g\rng \d\lambda(x)
\end{equation}
belongs to $\clm$ with $\nu^1(0)=\lambda$. Hence there exists a solution
to the martingale problem for $(A,\lambda)$ {\rm (}whose finite dimensional
distributions are $\{\nu^m\}${\rm )}.
\end{lemma}

\begin{proof}
It is easy to see that $\{\nu^m\}$ satisfy (\ref{r26}) since each
$\{\mu_x^m\}$ satisfies the same. Further, condition (\ref{r43}) on
$\lambda$ alongwith the definition of $\Phi^*_N$ implies that $\nu^1$
satisfies (\ref{r25}) and hence $\{\nu^m\}$ belongs to $\clm$. Thus the
corresponding process $Y$ is a solution to the martingale problem
for $(A,\lambda)$.
\qed
\end{proof}

We need one more observation on martingale problems before we can
state our result on $(T_t)$ defined by (\ref{r41}). \vspace{-.5pc}

\begin{lemma}\label{L8}
Let a process $X$ defined on $(\Omega,\clf,\IP)$ be a solution to
the $(A,\lambda)$ martingale problem and let $g$ be a
$[0,M]$-valued measurable function on $E$ {\rm (}where
$M<\infty${\rm )} such that $\lng \lambda ,g\rng=1$. Let $\gamma$
be defined by ${\d\gamma}/{\d\lambda}=g$. Let $\IQ$ be defined by
\begin{equation*}
\frac{\d\!\IQ}{\d\IP}=g(X_0).
\end{equation*}
Then{\rm ,} considered as a process on $(\Omega,\clf,\IQ\!\!)${\rm
,} $X$ is a solution to the $(A,\gamma)$ martingale problem.
\end{lemma}\vspace{-1pc}

\begin{proof}
Since $g$ is bounded it follows that
\begin{align*}
\IE_{\IQ}\left[\int_0^N\Phi(X_u)\d u\right]&\leq M\IE_{\IP}\left[\int_0^N\Phi(X_u)\d u\right]\\
&< \infty.
\end{align*}
Moreover, since ${\d\IQ}/{\d\IP}$ is $\sigma({X_0})$ measurable,
it follows that $f(X_t)-\int_0^tAf(X_u)\d u$ is a martingale on
$(\Omega,\clf,\IQ\!)$ (as it is a martingale on
$(\Omega,\clf,\IP)$). The result follows upon noting that
$\IQ\circ (X_0)^{-1}=\gamma $. \qed
\end{proof}

In addition to Assumption A3, we need to assume the following in
order to show that $\{T_t\}$ is a semigroup. \pagebreak

\bhyp\label{h7} {\rm There exists a sequence $\{h_n\hbox{:}\ n\geq
1\}$ of $[0,\infty)$-valued Borel measurable functions on $E$ such
that for every $\lambda\in\clp(E)$ satisfying\vspace{-.2pc}
\begin{equation}\label{r51}
\lng h_n,\lambda \rng < \infty\quad\forall n\geq 1,
\end{equation}
any two solutions to the $(A,\lambda)$ martingale problem that are
continuous in probability have the same finite-dimensional
distributions.}
\ehyp

Thus, in order to verify that Assumption A6 holds in a given example, we can show
that the uniqueness holds under finitely many (or even countably many)
integrability condition(s). We are now in a position to prove the
semigroup property of $(T_t)$. In the course of the proof, we also get,
with little extra work, the result that every solution to the martingale
problem satisfies the Markov property. The Markov property can also be
obtained by following arguments as in \cite{K1}.

\begin{theorem}[\!]\label{main}
Suppose that $A,\cld(A)$ satisfy Assumptions A{\rm 1,} A{\rm 2,}
A{\rm 3} and A{\rm 6}.
\begin{enumerate}
\leftskip .2pc
\renewcommand{\labelenumi}{\rm (\roman{enumi})}
\item The martingale problem for $(A,\lambda)$ is well-posed in the
class of solutions that are continuous in probability if and only
if $\lambda\in\clp_\Phi$. Further{\rm ,} the finite-dimensional
laws of the solution $Y$ that is continuous in probability are
given by {\rm (\ref{r44})}.

\item Let $\lambda\in\clp_\Phi$. Let $X$ be a solution to the
$(A,\lambda)$ martingale problem {\rm (}defined on some
probability space $(\Omega , \clf , \IP)${\rm )}. Further{\rm ,}
let $X$ be continuous in probability. Then $X$ is a Markov process
and the associated semigroup $\{T_t\hbox{\rm :}\ t\geq 0\}$ is
defined by {\rm (\ref{r41})}.\vspace{-.7pc}
\end{enumerate}
\end{theorem}

\begin{proof}$\left.\right.$%\vspace{.5pc}
\begin{enumerate}
\leftskip .2pc
\renewcommand{\labelenumi}{(\roman{enumi})}
\item Let $\lambda\in\clp_\Phi$. We have seen in Lemma \ref{L7} that the
$(A,\lambda)$ martingale problem admits a solution $X$ whose
finite-dimensional distributions are given by (\ref{r44}). Let $X$
be defined on $(\Omega , \clf , \IP)$. This process $X$ is
continuous in probability. Let $Y$ be another solution to the
$(A,\lambda)$ martingale problem defined on $(\tilde{\Omega} ,
\tilde{\clf} , \tilde{\IP})$ such that $Y$ is continuous in
probability. Define $g$ on $E$ by\vspace{-.2pc}
\begin{equation*}
\hskip -1.25pc g(x)=C\sum_{n=1}^\infty 2^{-n}\frac{1}{1+h_n(x)},
\end{equation*}
where $C$ is a constant that is chosen so that $\lng
\lambda,g\rng=1$. Define probability measures $\gamma$, $\IQ$ and
$\tilde{\IQ}$ by\vspace{-.2pc}
\begin{equation*}
\hskip -1.25pc
\frac{\d\gamma}{\d\lambda}=g,\;\;\frac{\d\IQ}{\d\IP}=g(X_0)\;\mbox
{ and } \; \frac{\d\tilde{\IQ}}{\d\tilde{\IP}}=g(Y_0).
\end{equation*}
\looseness -1 By Lemma \ref{L8}, $X$ on $(\Omega , \clf , \IQ)$
and $Y$ on $(\tilde{\Omega}, \tilde{\clf} , \tilde{\IQ})$ are
solutions to the $(A,\gamma)$ \hbox{martingale} problem. Further,
these processes are continuous in probability. By construction,
$\gamma$ satisfies (\ref{r51}) and hence by Assumption A6, the
finite-dimensional distributions of $X$ on \hbox{$(\Omega, \clf , \IQ)$}
are the same as those of $Y$ on \hbox{$(\tilde{\Omega}, \tilde{\clf},
\tilde{\IQ})$}. This~in~turn implies that the finite-dimensional
distributions of $X$ on \hbox{$(\Omega, \clf, \IP)$} are the same as those
of $Y$ on \hbox{$(\tilde{\Omega}, \tilde{\clf} , \tilde{\IP})$}.
This proves well-posedness of the martingale problem for
\hbox{$(A,\lambda)$}.\vspace{-.2pc}

Conversely, let $X$ be a solution of the $(A,\lambda)$ martingale
problem that is continuous in probability. This time define\vspace{-.2pc}
\begin{equation*}
\hskip -1.25pc g(x)=C\sum_{n=1}^\infty
2^{-n}\frac{1}{1+\Phi^*_n(x)},
\end{equation*}
$\left.\right.$\vspace{-1pc}

\pagebreak

\noindent where $C$ is a constant that is chosen so that $\lng \lambda,g\rng=1$.
Define probability measures $\gamma$ and $\IQ$ by
\begin{equation*}
\hskip -1.25pc \frac{\d\gamma}{\d\lambda}=g,\;\mbox { and
}\;\frac{\d\IQ}{\d\IP}=g(X_0).
\end{equation*}
By Lemma \ref{L8}, $X$ is a solution to the $(A,\gamma)$
martingale problem under $\IQ$ and $X$ is continuous in $\IQ$
probability. By the first part, we have that the regular
conditional probability distribution of $(X_{t_1},X_{t_2},\dots
,X_{t_m})$ given by $\sigma(X_0)$ is $\mu_{X_0}^m(t_1,t_2,\dots
,t_m)$. As a consequence
\begin{equation}\label{r61}
\hskip -1.25pc \IE_{\IQ}\left[\int_0^N\Phi(X_s)\d
s|\sigma(X_0)\right]=\Phi^*_N(X_0).
\end{equation}
Since ${\d\IQ}\!/{\d\IP}$ is $\sigma(X_0)$ measurable, (\ref{r61})
implies that
\begin{equation*}
\hskip -1.25pc \IE_\IP\left[\int_0^N\Phi(X_s)\d
s|\sigma(X_0)\right]=\Phi^*_N(X_0)
\end{equation*}
and hence
\begin{align*}
\hskip -1.25pc \IE_\IP\left[\int_0^N\Phi(X_s)\d s\right]&=\IE_\IP\left[\Phi^*_N(X_0)\right]\\
\hskip -1.25pc &=\lng \Phi^*_N, \lambda\rng.
\end{align*}
Since $X$ is a solution to the $(A,\lambda)$ martingale problem, the LHS
above is finite for all $N$ and hence $\lambda\in\clp_\Phi$.

\item Let $X$ be a solution to the $(A,\lambda)$ martingale problem that
is continuous in probability (defined on some probability space
$(\Omega, \clf , \IP)$). Fix $m\geq 1$ and $0\leq u_1<u_2< \dots <
u_m\leq s$ and $h_1,h_2\dots h_m$ bounded positive continuous
functions. Define a probabiltiy measure $\IQ$ on $(\Omega , \clf ,
\IP)$ by
\begin{equation*}
\hskip -1.25pc \frac{\d\IQ}{\d\IP}=Ch_1(X_{u_1})h_2(X_{u_2})\dots
h_m(X_{u_m}),
\end{equation*}
where the constant $C$ is chosen such that $\IQ$ is a probability measure. Define $Y$ by
\begin{equation*}
\hskip -1.25pc Y_t=X_{s+t},\;\;\;t\geq 0.
\end{equation*}
Then using ${\d\IQ}/{\d\IP}$ which is bounded (say by $M$), we get
\begin{align}
\hskip -1.25pc \IE_{\IQ}\left\lbrack \int_0^T \Phi (Y_u)\d u\right\rbrack &=  \IE_{\IQ}\left\lbrack\int_s^{T+s} \Phi (X_u)\d u\right\rbrack\nonumber\\[.3pc]
\hskip -1.25pc &\leq M\IE_\IP\left\lbrack\int_s^{T+s} \Phi (X_u)\d
u\right\rbrack <\infty.
\end{align}
Further, it can be shown that $Y$ is a solution to the
$(A,\gamma)$ martingale problem where $\gamma = \IQ \circ
[Y(0)]^{-1}$. Of course, $Y$ is continuous in probability. Hence,
by part (i) above we get that $\gamma\in \clp_\Phi$ and that the
finite-dimensional distributions are given by (\ref{r44}) (with
$\lambda$ replaced by $\gamma$). Thus, for $g_1,\dots, g_k\in
C_b(E)$ and $0\leq s_1 < \dots <s_k$,
\begin{align*}
\hskip -1.25pc \IE_{\IQ}\left[g_1(Y_{s_1})\dots g_k(Y_{s_k})\right]&=\int\lng\mu_x^k(s_1,\dots ,s_k),g_1\otimes\cdots\otimes g_k\rng \d\gamma(x)\\[.3pc]
\hskip -1.25pc &=\IE_{\IQ}\big[\lng\mu_{Y_0}^k(s_1,\dots
,s_k),g_1\otimes\cdots\otimes g_k\rng \big]
\end{align*}
and so (using $k=1$, $s_1=t$ and $g_1=g$) we can conclude that
 \begin{align*}\label{r73}
\hskip -1.25pc &\IE_{\IP}\left[Ch_1(X_{u_1})h_2(X_{u_2})\ldots h_m(X_{u_m})g(X_{s+t})\right]\\[.3pc]
\hskip -1.25pc &\quad\,=\IE_{\IQ}\left[g(Y_t)\right] = \IE_{\IQ}\big[\lng\mu_{Y_0}^1(t),g\rng\big]\\[.3pc]
\hskip -1.25pc
&\quad\,=\IE_{\IP}\big[Ch_1(X_{u_1})h_2(X_{u_2})\ldots
h_m(X_{u_m})\lng\mu_{X_s}^1(t),g\rng\big]
\end{align*}
for all $0\leq u_1<u_2< \cdots < u_m\leq s$ and $h_1,h_2\ldots
h_m$ bounded positive continuous functions, $m\geq 1$. As a
consequence,
\begin{equation*}
\hskip -1.25pc
\IE_{\IP}\left[g(X_{s+t})\big\vert\sigma(X_u\hbox{:}\ 0\leq u\leq
s)\right]=\lng\mu_{X_s}^1(t),g\rng=(T_tg)(X_s).
\end{equation*}
This completes the proof.\qed\vspace{-2pc}
\end{enumerate}
\end{proof}

\section{Criterion for an invariant measure}

Several papers gave criterion for a measure to be invariant for the
semigroup $(T_t)$ arising from a well-posed martingale problem
\cite{BB,RLKAB1,RLKAB2,RLKAB3,Ech,KS1,KS2}. These papers assumed
different sets of conditions on $(A, D(A))$. It was shown that existence
of solution for each degenerate initial and
\begin{equation*}
\int (Af)\d\lambda=0\quad\forall f\in \cld(A)
\end{equation*}
gives existence of a stationary solution of the martingale problem for
$(A,\lambda)$. In addition, if the martingale problem is well-posed and
there is a semigroup $(T_t)$ associated with it, it follows that
$\lambda$ is an invariant measure for $(T_t)$.

Well-posedness of the martingale problem in the class of r.c.l.l.
solutions is sufficient for the existence of the semigroup $(T_t)$ (see
Theorem 4.4.6 of \cite{EK}).

In the light of the results obtained in the previous section, we can improve
on this criterion for invariant measure.

We introduce another condition on $A$ and $\Phi$ (appearing in
Assumption A1).

\bhyp
{\rm $\Phi $ and  $Af$,  for every  $f\in\cld(A)$, are continuous.}
\ehyp

\setcounter{defin}{0}
\begin{lemma}\label{L9}
Suppose that $A,\cld(A)$ satisfy Assumptions A{\rm 1,} A{\rm 2,}
A{\rm 3} and A{\rm 7}. Then $A$ satisfies the positive maximum
principle{\rm ,} i.e. if $f\in\cld(A)$ and $z\in E$ are such that
$f(z)\geq 0$ and $f(z)\geq f(x)$ for all $x\in E${\rm ,} then
\begin{equation*}
Af(z)\leq 0.
\end{equation*}
\end{lemma}

\begin{proof}
Let $X$ be a solution to $(A,\delta_z)$ martingale problem defined on
$(\Omega,\clf,\IP)$ that is continuous in probability. Let
$\clf_t={}^*\clf^X_t$ and
\begin{equation*}
M_t=f(X_t)-\int_0^tAf(X_u)\d u.
\end{equation*}
Then $(M_t,\clf_t)$ is a martingale. Let $\sigma_t$, $0\leq
t<\infty$ be the increasing family of $(\clf_t)$ stopping times
defined by
\begin{equation*}
\sigma_t=\inf\left\lbrace s\geq 0\hbox{:}\
\int_0^s(1+\Phi(X_u))\d u\geq t\right\rbrace.
\end{equation*}
Note that $\sigma_t\leq t$ for all $t$. Since
$\IE_\IP[\int_0^s\Phi(X_u)\d u]<\infty$, it follows that $\sigma_t$
increases to $\infty$ a.s..

Let $N_t=M_{\sigma_t}$, $Y_t=X_{\sigma_t}$ and $\clg_t=\clf_{\sigma_t}$.
Then, it follows that $(N_t,\clg_t)$ is a local martingale. Moreover,
$t\longrightarrow \sigma_t$ is continuous and hence $Y$ is also
continuous in probability. Using change of variable, it is easy to see
that
\begin{equation*}
N_t=f(Y_t)-\int_0^t\frac{Af(Y_r)}{1+\Phi(Y_r)}\d r.
\end{equation*}
Since $Af(x)\leq C_f\Phi(x)$, it follows that $N$ is bounded and
hence is a martingale. Since $f$ has a maximum at $z$ and
\begin{equation*}
\IE_\IP\left[f(Y_t)-f(z)-\int_0^t\frac{Af(Y_r)}{1+\Phi(Y_r)}\d r\right]=0,
\end{equation*}
it follows that (using Fubini's theorem)
\setcounter{equation}{0}
\begin{equation}\label{r77}
\int_0^t\IE_\IP\left[\frac{Af(Y_r)}{1+\Phi(Y_r)}\right]\d r\leq
0\quad\forall t>0.
\end{equation}
Since $Y$ is continuous in probability and $Af(x)\leq C_f\Phi(x)$,
it follows that
\begin{equation*}
r\longrightarrow \IE_\IP\left[\frac{Af(Y_r)}{1+\Phi(Y_r)}\right]
\end{equation*}
is continuous. Now dividing the LHS in (\ref{r77}) by $t$ and taking limit as $t\raro 0$ we get
\begin{equation*}
\frac{Af(z)}{1+\Phi(z)}\leq 0.
\end{equation*}
Since $\Phi(z)\geq 0$ this completes the proof. \qed

\end{proof}

Here is yet another assumption on $A,\cld(A).$

\bhyp
{\rm $\cld(A)$ is an algebra that contains constants and separates points in $E$.}
\ehyp

\begin{theorem}[\!]
Suppose that $A,\cld(A)$ satisfy Assumptions A{\rm 1,} A{\rm 2,}
A{\rm 3,} A{\rm 6,} A{\rm 7} and A{\rm 8}. Let $(T_t)$ be the
semigroup associated with $(A,\cld(A))$ by Theorem {\rm \ref{main}}.

If $\lambda\in\clp(E)$ is such that $\int \Phi \d\lambda<\infty$ and
\begin{equation}\label{inv}
\int (Af)(x)\d\lambda(x)=0\quad\forall f\in\cld(A),
\end{equation}
then $\lambda$ is an invariant measure for the semigroup $(T_t)$ and the
solution to the $(A,\lambda)$ martingale problem that is continuous in
probability is a stationary process.
\end{theorem}

\begin{proof}
In view of Lemma \ref{L9} and the assumptions made in the statement of
this theorem, the proof of Theorem 3.1 in \cite{RLKAB2} gives the existence of a stationary solution to the $(A,\lambda)$
martingale problem. Since the solution (say $X$) is stationary, the
mapping $t\longrightarrow \cll(X_t)$ is continuous (it is a constant)
and hence by Theorem \ref{NT1}, $X$ is continuous in probability. Now,
Theorem \ref{main} implies that $\lambda$ is an invariant measure for
$(T_t)$.
\qed
\end{proof}

\begin{remark}
{\rm The criterion for invariant measure given above is true even if
Assumption A7 above is not true but instead one assumes that the
operator $A$ satisfies Assumptions A9, A10 and A11 given below. This is helpful,
e.g., when $Af$ is allowed to be a discontinuous function (see
\cite{RLKAB3,KS1}).

\bhyp
{\rm $A$ satisfies the positive maximum principle.}
\ehyp

\bhyp {\rm There exists a complete separable metric space $U$, an
operator $\hat{A}\hbox{:}\  \cld (A)\raro C(E\times U)$ and a
transition function $\eta$ from $(E,\clb(E))$ into $(U,\clb(U))$
such that
\begin{equation}\label{b1}
(Af)(x) = \int_U \hat{A}f(x,u)\eta(x,\d u).
\end{equation}}
\ehyp

\bhyp
{\rm There exists $\hat{\Phi} \in C(E\times U)$ such that for all
$f\in\cld (A)$, there exists $C_f<\ity$ satisfying
\begin{align}
|\hat{A}f(x,u)| &\leq C_f\hat{\Phi}(x,u) \quad \fal x,u\in E\times U,\label{rk0}\\[.5pc]
\Phi (x) &= \int_U \hat{\Phi} (x,u)\eta(x,\d u) < \ity.\label{b2}
\end{align}}
\ehyp

Under these conditions, existence of a stationary solution to the
$(A,\lambda)$ martingale problem was proven in \cite{RLKAB3}. Rest of
the argument is as in the proof of the above theorem.}
\end{remark}

\section{Example}

We give two examples of processes that are continuous in probability and
which arise as solutions of well-posed martingale problems but such that
they do not admit any r.c.l.l. modification. The results of the previous
section, however, are applicable.

\begin{exa}
{\rm Let $ E = [0,1)$. Let $\cld (A)$ be the class of functions $f$
that are restrictions of some periodic function $g \in C^2_b(\IR)$ with
period $1$. Further for $f \in \cld (A)$ define $Af$ by $Af =
\frac{1}{2} f^{\prime\prime}$. Then $A$ and $\cld(A)$ satisfy the

conditions of Theorems \ref{mainT} and \ref{main}.

It follows easily that if $W$ is a one-dimensional standard
Brownian motion then $X_t = W_{t}\! \pmod 1$ is a solution to the
martingale problem for $A$. Moreover, for any other solution~$Y$
of the martingale problem, it is easy to check that $Y$ behaves
like a Brownian motion as long as it does not hit the boundary.
Now, uniqueness can be shown using localisation arguments as in
Theorem 6.6.1 of \cite{SV}.

Note that almost every path of the unique solution $X$ is neither
r.c.l.l. nor l.c.r.l. However the set of discontinuity points of $X$ is
contained in the set
\begin{equation*}
\{ t\hbox{:}\ W_t \mbox{ is an integer}\}.
\end{equation*}
This implies that $X$ is continuous in probability.}
\end{exa}

\begin{exa}
{\rm Let $E = (0,\infty )$ and let $\mu$ be a probability measure on $E$ with
$\mu \{(0, a)\} = 0$ for some $a>0$. Let $\cld (A)$ be defined by
\begin{equation*}
\cld (A) = \left\{ f\in C^2_b(E)\hbox{:}\  \lim_{x\rightarrow 0}
f(x) = \int f\;\d\mu \right\}.
\end{equation*}
For $f\in \cld(A)$ define $Af$ by $Af = \frac{1}{2}f^{\prime\prime}$.
Once again, $A$ and $\cld(A)$ satisfy the conditions of Theorems \ref{mainT} and \ref{main}.

The uniqueness of solution for the martingale problem for $A$ can
also be shown using localisation arguments as in Theorem 6.6.1 of
\cite{SV}. To construct the unique solution for the
$(A,\delta_x)$ martingale problem we can proceed as follows.

Let $\{W^{z,i}\hbox{:}\ i\geq 0\}$ be independent one-dimensional
standard Brownian motions starting at $z$. Define
\begin{equation*}
\tau^{z,i} = \inf \{t>0\hbox{:}\ W^{z,i} = 0\}.
\end{equation*}
Note that
\setcounter{equation}{0}
\begin{equation}\label{tauzi}
\tau^{z,i} < \infty \mbox { a.s. for every } z, i.
\end{equation}
Let $U_1, U_2, \dots$ be i.i.d. random variables with common
distribution $\mu$ and which are independent of all
$\{W^{z,i}\hbox{:}\ i\geq 1\}$. Define
\begin{equation*}
X^x_t = \begin{cases}
           W^{x,0}_t,  & \text{for $t < \tau^{x,0}$}, \\[.5pc]
           W^{U_i,i}_t,  & \text{for $\tau^{U_{i-1},i-1} \leq t < \tau^{U_i,i}; \quad i\geq 1$}.
           \end{cases}
\end{equation*}
Then it is easily checked that $X^x$ is a solution of the martingale
problem for $A$ starting at $x$ and which is also continuous in
probability. Then \eqref{tauzi} and the fact that $0\not\in E$ together
imply that almost every path of $X^x$ is not left continuous. }
\end{exa}

\end{document}